\documentclass[titlepage,11pt]{article}
\usepackage{latexsym}
\usepackage{amsfonts}
\usepackage{amsmath}
\usepackage[T1]{fontenc}
\usepackage{mathtools}
\newcommand\blackslug{\hbox{\hskip 1pt \vrule width 4pt height 8pt depth 1.5pt
        \hskip 1pt}}
\newcommand\bbox{\hfill \quad \blackslug \bigbreak}

\def\LL{,\ldots,}
\def\cupcup{\cup\cdots\cup}

\title{Polynomial bounds for chromatic number.\\ V. Excluding a tree of radius two and a complete multipartite graph}
\author{Alex Scott\thanks{Research supported by EPSRC grant EP/V007327/1.}\\
Mathematical Institute, University of Oxford, Oxford OX2 6GG, UK
\\
\\
Paul Seymour\thanks{Supported by AFOSR grant
A9550-19-1-0187, and by NSF grant  DMS-1800053.}\\
Princeton University, Princeton, NJ 08544}

\date{}

\newtheorem{thm}{}[section]

\newcommand{\Proof}{\noindent{\bf Proof.}\ \ }

\begin{document}
\maketitle
\begin{abstract}
The Gy\'arf\'as-Sumner conjecture says that for every forest $H$ and every integer $k$, if $G$ is $H$-free and does not contain 
a clique on $k$ vertices then it has bounded chromatic number.  (A graph is {\em $H$-free} if it does not contain an induced copy 
of $H$.) 
Kierstead and Penrice proved it for trees of radius at most two, but otherwise the conjecture is known only for a few simple
types of forest.  More is known if we exclude a complete bipartite 
subgraph instead of a clique: R\"odl showed that, for every forest $H$, if $G$ is $H$-free and does not contain $K_{t,t}$ as a subgraph then it has 
bounded chromatic number.  In an earlier paper with Sophie Spirkl, we strengthened R\"odl's result, showing that for every forest $H$, 
the bound on chromatic number can be taken to be polynomial in $t$.  In this paper, we prove a related strengthening of the 
Kierstead-Penrice theorem, showing that for every tree $H$ of radius two and integer $d\ge 2$, if $G$ is $H$-free and does not 
contain as a subgraph the complete $d$-partite graph with parts of cardinality $t$, then its chromatic number is at most polynomial in $t$.

\end{abstract}
\section{Introduction}

The Gy\'arf\'as-Sumner conjecture~\cite{gyarfas,sumner} says:
\begin{thm}\label{GSconj}
{\bf Conjecture: }For every forest $H$ there is a function $f$ such that $\chi(G)\le f(\omega(G))$ for every $H$-free graph $G$.
\end{thm}
($G$ is {\em $H$-free} if no induced subgraph of $G$ is isomorphic
to $H$; and $\chi(G),\omega(G)$ denote the chromatic number and the size of the largest clique of $G$, respectively.)

This is open in general, although it is known to hold~\cite{scott} for graphs that do not contain any induced {\em subdivision} of 
$H$, and has been proved for a few special kinds of forest.  Notably, Kierstead and Penrice~\cite{kierstead} proved:
\begin{thm}\label{kierstead}
	For every tree $H$ of radius two, there is a function $f$ such that $\chi(G)\le f(\omega(G))$ for every $H$-free graph $G$.
\end{thm}

These statements can also be phrased in terms of $\chi$-bounded classes.  A class of graphs is {\em hereditary} if it is closed 
under taking induced subgraphs; and a hereditary class $\mathcal G$ of graphs is {\em $\chi$-bounded} if there is a function $f$ 
such that $\chi(G) \le f(\omega(G))$ for every graph $G\in\mathcal G$.  Thus conjecture \ref{GSconj} says that, for every forest 
$H$, the class of $H$-free graphs is $\chi$-bounded; and \ref{kierstead} says that the class of $H$-free graphs is $\chi$-bounded 
when $H$ is a tree of radius two.

There has been a great deal of recent progress on $\chi$-bounded classes (see~\cite{survey} for a survey).  In most cases, the 
proofs give bounds on the chromatic number that grow relatively quickly (often superexponentially) in the clique number.  
However, a striking conjecture of Esperet~\cite{esperet} asserts that this is not necessary, and that for every $\chi$-bounded 
class, the function $f$ can be taken to be polynomial.  Esperet's conjecture has been shown to be false in its full generality~\cite{brianski}, 
but remains open for classes of graphs excluding a forest; and in that case, the Gy\'arf\'as-Sumner conjecture and 
Esperet's conjecture would together give the following:
\begin{thm}\label{Econj}
{\bf Conjecture: }For every forest $H$, there is a polynomial $f$ such that $\chi(G)\le f(\omega(G))$ for every $H$-free graph $G$.
\end{thm}
So far, this is only known for a few classes of trees (see \cite{poly6, liu,poly2,poly3,Schiermeyer}).

While the conjectures \ref{GSconj} and \ref{Econj} remain open, more is known if we exclude a complete bipartite graph rather than a clique.  
R\"odl (see~\cite{gst,kr}) proved that:

\begin{thm}\label{rodl}
For every forest $H$ and integer $t\ge2$, there exists $k$ such that if $G$ is $H$-free and does not contain $K_{t,t}$ as 
a subgraph then  $\chi(G)\le k$.
\end{thm} 

It will be helpful to define one piece of notation.
For a graph $G$, and integer $d\ge 1$, let $\tau_d(G)$ denote the largest $t$ such that $G$ has a subgraph (not necessarily induced)
isomorphic to the complete $d$-partite graph with each part of cardinality $t$.  Thus $\tau_1(G)=|G|$, and $\tau_2(G)$ is the 
largest $t$ such that $G$ contains $K_{t,t}$ as a subgraph, and \ref{rodl} says that for every forest $H$ there is a 
function $f$ such that every $H$-free graph $G$ satisfies $\chi(G)\le f(\tau_2(G))$.
It is natural to ask whether $f$ can be taken to be a polynomial in R\"odl's result.  When $H$ is a path, this was proved by 
Bonamy, Bousquet, Pilipczuk, Rz\k{a}\.zewski, Thomass\'e and Walczak~\cite{bonetc}.   We proved the general case with Sophie Spirkl in~\cite{poly1}:
\begin{thm}\label{poly1}
For every forest $H$, there exists $c>0$ such that $\chi(G)\le \tau_2(G)^c$ for every $H$-free graph $G$.
\end{thm} 
Note that this is a special case of \ref{Econj}, as $\tau_2(G)\ge \lfloor\omega(G)/2\rfloor$.  \ref{Econj} would also imply 
that the same is true for $\tau_d(G)$ instead of $\tau_2(G)$ for any fixed value of $d\ge 2$ (except if $\tau_d(G)\le 1$)
since $\tau_d(G)\ge \lfloor \omega(G)/d\rfloor$.
This has not been proved in general~--~indeed, proving it for a forest $H$ would show that $H$ also satisfies 
the Gy\'arf\'as-Sumner conjecture~--~but in this paper
we prove it when $H$ is a tree of radius two.
Our main result is the following extension of \ref{kierstead}:

\begin{thm}\label{mainthm}
        For every tree $H$ of radius two, and every integer $d\ge 1$, there is a polynomial $f$ such that 
	$\chi(G)\le f(\tau_d(G))$ 
	for every $H$-free graph $G$.
\end{thm}

A referee suggests two further open questions on these lines:
\begin{itemize}
\item (Extending \ref{mainthm} to other trees $H$ that we know satisfy \ref{GSconj}.) Is it true that if $H$ is a path, then for every integer $d\ge 1$, there is a polynomial $f$ such that
        $\chi(G)\le f(\tau_d(G))$
        for every $H$-free graph $G$?
\item (An analogue of Esperet's (false) conjecture for $\tau_d(G)$.) Let $\mathcal{C}$ be a hereditary class of graphs, 
and $d\ge 1$. Suppose that there is a function $f$
such that $\chi(G)\le f(\tau_d(G))$ for each $G\in \mathcal{C}$. Can we always choose $f$ to be a polynomial? What if $d=2$?
\end{itemize}
We note that for the five-vertex path $P_5$, the best current upper bound on chromatic number is $\omega^{\log_2 \omega}$ \cite{poly4}.
A polynomial upper bound, as asserted by conjecture \ref{Econj}, would imply that $P_5$ satisfies the
Erd\H{o}s-Hajnal conjecture~\cite{EH0, EH} ($P_5$ is currently the smallest open case of the 
Erd\H{o}s-Hajnal conjecture, after $C_5$ was recently proved in~\cite{fivehole}).  Since $P_5$ is a tree of radius two, it would be very nice if the function $f$
in \ref{mainthm} had polynomial dependence on $d$. But the function $f$ we prove in this paper has doubly-exponential
dependence on $d$. Incidentally, if we take $t=1$ in \ref{mainthm}, we have proved that $\chi(G)\le f(\omega(G))$
for every $H$-free graph $G$, where $f$ is doubly-exponential in $\omega(G)$. While this is admittedly 
fast-growing,
the bound is much smaller than that of Kierstead and Penrice~\cite{kierstead}.  

We use standard notation.  For a graph $G$, we denote the number of vertices 
by $|G|$.
When $X\subseteq V(G)$, $G[X]$ denotes the subgraph induced on $X$. We write $\chi(X)$ for $\chi(G[X])$
when there is no ambiguity. If $v\in V(G)$, a {\em non-neighbour} of $v$ in $G$ means a vertex $u$ of $G$ different from
 $v$ and nonadjacent to $v$.

\section{Some Ramsey-type lemmas}

We will use the following well-known version of Ramsey's theorem, proved (for instance)
in~\cite{poly2}:
\begin{thm}\label{ramsey}
        Let $x\ge 2$ and $y\ge 1$ be integers. For a graph $G$, if $|G|\ge x^y$, then $G$ has either a clique of
        cardinality $x+1$, or a stable set of cardinality $y$.
\end{thm}

We also need the next result:
\begin{thm}\label{levels2}
        Let $s\ge 2$ and $t\ge 1$ be integers, and
        let $G$ be a graph with $\tau_{d+1}(G)<t$.
        Let $L_1\LL L_{s^{2d+2}}$ be pairwise disjoint subsets of $V(G)$, each of
        cardinality at least $2^{s^{2d+2}}t^{ds+s^2+s}$. Then there exist $I\subseteq \{1\LL s^{2d+2}\}$ with $|I|=s$,
        and a subset $X_i\subseteq L_i$
        for each $i\in I$, where
        $\bigcup_{i\in I}X_i$ is a stable set, and $|X_i|\ge s$ for each $i\in I$.
\end{thm}

For inductive purposes, we will prove the following stronger (but messier) form: \ref{levels2} follows from it 
by substituting $a=b=c=s$.
\begin{thm}\label{levels}
	Let $a\ge 2$ and $t\ge 1$ be integers. For all integers $b,c,d\ge 0$ with $b\le a$ and $c\ge 1$, define
	$$k_{b,c,d}=\begin{cases*}
		(ac)^{d+1} & \text{ if } $b=a$\\ 
		b(ac)^{d}+(a(c-1))^{d+1}+1 & \text{ if $b<a$ and $d>0$}\\
		1 & \text{ if $b<a$ and  $d=0$}
	\end{cases*}$$
	Define $p_{b,c,d} = 2^{k_{b,c,d}}t^{a(c+d)+b}$. 
	Now let $b,c,d\ge 0$ be integers with $b\le a$, and let $G$ be a graph with $\tau_{d+1}(G)<t$. 
	Let $L_1\LL L_{k_{b,c,d}}$ be pairwise disjoint subsets of $V(G)$, each of 
	cardinality at least $p_{b,c,d}$. Then there exist $I\subseteq \{1\LL k_{b,c,d}\}$ with $|I|=c$, 
	and a stable subset $X\subseteq L_1\cupcup L_{k_{b,c,d}}$, where
        $|X\cap L_i|\ge a$ for each $i\in I\setminus \{1\}$,
	and $|X\cap L_1|\ge b$ if $1\in I$.
\end{thm}
\Proof
We proceed by induction on $(a+1)(c+d)+b$ (the numbers $a,t$ are fixed throughout the proof).
If $d=0$ then $\tau_{d+1}(G)=|G|$, so $|G|<t\le p_{b,c,d}$, and there is no choice of $L_1$ satisfying the hypothesis, 
and therefore the theorem holds. So we may assume that $d\ge 1$.

Suppose that $c=1$. Since $G$ has no clique of cardinality $(d+1)t$ (because $t>\tau_{d+1}(G)$), \ref{ramsey} implies that
every set of $((d+1)t)^b$ vertices of $G$ includes a stable set of cardinality $b$. 
Thus it suffices to show that $p_{b,c,d}\ge ((d+1)t)^b$ when $c=1$, that is, we must show that
$$2^{k_{b,1,d}}t^{a(1+d)+b}\ge ((d+1)t)^b.$$
Since $t^{a(1+d)+b}\ge t^b$, it is enough to show that 
$2^{k_{b,1,d}}\ge (d+1)^b$. But $k_{b,1,d}\ge ba^d$, so it suffices to show that 
$2^{ba^d}\ge (d+1)^b$, that is,
$a^d \ge \log_2(d+1)$. Since $a\ge 2$, and $2^d\ge  \log_2(d+1)$, this is true, so we may assume that $c\ge 2$.

Suppose that $b=0$, and therefore $k_{b,c,d} = k_{a,c-1,d}+1$.
By applying the inductive hypothesis to $L_2\LL L_{k_{b,c,d}}$, with $b,c,d$ replaced by
$a,c-1,d$ respectively,
we deduce that there exist $I'\subseteq \{2\LL k_{b,c,d}\}$ with $|I'|=c-1$, and a stable subset $X\subseteq L_2\cupcup L_{k_{b,c,d}}$,
where $|X\cap L_i|\ge a$ for each $i\in I'$.
Then setting $I=I'\cup \{1\}$ satisfies the theorem. Thus we may assume that $b\ge 1$.
\\
\\
(1) {\em The following inequalities hold:
\begin{align*}
	k_{b,c,d}- k_{b-1,c,d} &\ge k_{a,c,d-1} \\
	p_{b,c,d}&\ge  2(d+1)^2t^2\\
	p_{b,c,d} &\ge   2(d+1)tp_{b-1,c,d}\\
	p_{b,c,d}&\ge  2^{k_{b,c,d}} t\\
	  p_{b,c,d}&\ge  tp_{b-1,c,d}+ p_{a,c,d-1}.
\end{align*}
}
\noindent
The first is clear (and holds with equality) if $b<a$, so we assume that $b=a$. Since $$1-1/c> (1-1/c)^{d+1}$$ (because $c\ge 2$ and $d\ge 1$),
we have $$(ac)^{d+1}> a(ac)^{d}+(a(c-1))^{d+1},$$
and so 
$$k_{a,c,d}\ge  a(ac)^{d}+(a(c-1))^{d+1}+1= k_{a-1,c,d}+k_{a,c,d-1},$$ 
and the first inequality follows.

For the second, we must show that $2^{k_{b,c,d}}t^{a(c+d)+b}\ge 2(d+1)^2t^2$. 
Since $t^{a(c+d)+b}\ge t^2$, it suffices to show that $2^{k_{b,c,d}}\ge 2(d+1)^2$, and this is true since 
$k_{b,c,d}\ge (ac)^{d}+1\ge 2^{d}+1$, and $2^{2^{d}+1}\ge 2(d+1)^2$. 
This proves the second inequality.

For the third,  we must show that
$$2^{k_{b,c,d}}t^{a(c+d)+b}\ge 2(d+1)t 2^{k_{b-1,c,d}}t^{a(c+d)+b-1},$$
that is, 
$$k_{b,c,d}-k_{b-1,c,d}\ge 1+\log_2(d+1).$$
But (using the first inequality if $b=a$), 
$k_{b,c,d}-k_{b-1,c,d} \ge (ac)^{d}\ge 2^{d}\ge  1+\log_2(d+1)$ as required. 

For the fourth, we must show that $2^{k_{b,c,d}}t^{a(c+d)+b}\ge  2^{k_{b,c,d}} t$,
which is clear. Finally, for the fifth, we must show that
$$2^{k_{b,c,d}}t^{a(c+d)+b}\ge t 2^{k_{b-1,c,d}}t^{a(c+d)+b-1} + 2^{k_{b,c,d-1}}t^{a(c+d-1)+b},$$
that is,
$$2^{k_{b,c,d}}\ge  2^{k_{b-1,c,d}} + 2^{k_{b,c,d-1}}t^{-a}.$$
Since $t\ge 1$, it suffices to show that
$2^{k_{b,c,d}}\ge  2\cdot 2^{k_{b-1,c,d}}$
and $2^{k_{b,c,d}}\ge 2\cdot 2^{k_{b,c,d-1}}$, that is,
$k_{b,c,d}> k_{b-1,c,d}$ and $k_{b,c,d}> k_{b,c,d-1}$, which are both true (since $ac\ge 2$). This proves (1).

\bigskip

Choose a clique $Y\subseteq L_1$, maximal such that at most $|Y|p_{b,c,d}/(2(d+1)t)$ vertices in $L_1$ have a non-neighbour in $Y$. (Possibly $Y=\emptyset$.)
Let $N$ be the set of vertices in $L_1\setminus Y$
that are adjacent to every vertex in $Y$.
Then:
\\
\\
(2) {\em $|N|\ge p_{b,c,d}/2$, and every vertex $v\in N$ has more than $p_{b,c,d}/(2(d+1)t)$ non-neighbours in $N$.}
\\
\\
Since $t> \tau_{d+1}(G)$, it follows that $G$ has no clique of cardinality $(d+1)t$, and so $|Y|<(d+1)t$.
Let $M=L_1\setminus (N\cup Y)$. Thus $|M|\le |Y|p_{b,c,d}/(2(d+1)t)$ from the choice of $Y$, and so 
\begin{align*}
	|Y\cup M|&\le |Y|\left(1+\frac{p_{b,c,d}}{2(d+1)t}\right)\le \left((d+1)t-1\right)\left(1+\frac{p_{b,c,d}}{2(d+1)t}\right)\\
	&=(d+1)t-1+\frac{p_{b,c,d}}{2} - \frac{p_{b,c,d}}{2(d+1)t}\le \frac{p_{b,c,d}}{2}
\end{align*}
since $2(d+1)^2t^2\le p_{b,c,d}$ by (1). Consequently $|N|\ge p_{b,c,d}/2$. This proves the first assertion. For the second, if some
vertex $v\in N$ has at most $p_{b,c,d}/(2(d+1)t)$ non-neighbours in $N$, 
then adding $v$ to $Y$ gives a  set $Y'$ such that at most $|Y'|p_{b,c,d}/(2(d+1)t)$ vertices in $L_1$ have a non-neighbour in $Y'$,
contrary to the maximality of $Y$. This proves (2).

\bigskip

We may assume that 
\\
\\
(3) {\em For each $v\in N$, there are fewer than $k_{b-1,c,d}$ values of $i\in \{2\LL k_{b,c,d}\}$ such that $v$ has at least $p_{b-1,c,d}$
non-neighbours in $L_i$.}
\\
\\
Suppose that there exists $I'\subseteq \{2\LL k_{b,c,d}\}$ with $|I'|=k_{b-1,c,d}$, such that  for each $i\in I'$ there is a set $L_i'\subseteq L_i$
of non-neighbours of $v$ with $|L_i'|\ge p_{b-1,c,d}$. Let $L_1'$ be the set of non-neighbours of $v$
in $N$; then $|L_1'|\ge p_{b-1,c,d}$ by (2), since $p_{b,c,d}/(2(d+1)t) \ge  p_{b-1,c,d}$ by (1). 
From the inductive hypothesis applied
to $L_i'\;(i\in I'\cup \{1\})$, with $b,c,d$ replaced by $b-1,c,d$ respectively, it follows that 
there exist $I\subseteq I'\cup \{1\}$ with $|I|=c$, and a stable subset $X\subseteq \bigcup_{i\in I}L_i'$, where
$|X\cap L_i'|\ge a$ for each $i\in I\setminus \{1\}$, and
$|X\cap L_1'|\ge b-1$ if $1\in I$.
If $1\notin I$ then the theorem holds. If $1\in I$, then by adding $v$ to $X$ we see that again
the theorem holds. This proves (3).

\bigskip

For each $v\in N$, let $I_v$ be the set of values of $i\in \{2\LL k_{b,c,d}\}$ such that $v$ has fewer than $p_{b-1,c,d}$
non-neighbours in $L_i$. Thus $|I_v|\ge k_{b,c,d}- k_{b-1,c,d}$ for each $v$, by (3). Since there are at most  $2^{k_{b,c,d}-1}$ choices of the set $I_v$,
and $|N|\ge p_{b,c,d}/2\ge 2^{k_{b,c,d}-1}t$ by (1), there exists $T\subseteq N$ with $|T|= t$ such that the sets $I_v\;(v\in T)$ are all 
equal, and equal to some $I'$ say. For each $i\in I'$, since each $v\in T$ has at most $p_{b-1,c,d}$ non-neighbours in $L_i$,
it follows that there are at most $tp_{b-1,c,d}$  vertices in $L_i$ that have a non-neighbour in $T$. Since 
$|L_i|\ge p_{b,c,d}\ge tp_{b-1,c,d}+ p_{a,c,d-1}$ by (1),
there is a subset $L_i'\subseteq L_i$ with $|L_i'|\ge p_{a,c,d-1}$, such that every vertex in $T$ is adjacent to every vertex of $L_i'$.
Since $\tau_{d+1}(G)<t$, it follows that $\tau_{d}(G')<t$, where $G'=G[\bigcup_{i\in I'}L_i']$. 
Since $|I|\ge k_{b,c,d}- k_{b-1,c,d} \ge k_{a,c,d-1}$, the inductive hypothesis (on $c+d$) applied to $L_i'\;(i\in I')$, with $b,c,d$ replaced by 
$a,c,d-1$ respectively, implies that 
there exist $I\subseteq I'$ with $|I|=c$, and a stable subset $X\subseteq \bigcup_{i\in I} L_i'$, where
$|X\cap L_i'|\ge a$ for each $i\in I$.
This proves \ref{levels}.~\bbox

If $v\in V(G)$ and $B\subseteq V(G)$ with $v\notin B$, we say $v$ is {\em complete} to $B$ if $v$ is adjacent to 
every vertex in $B$, and $v$ is {\em anticomplete} to $B$ if $v$ has no neighbours in $B$. If $A,B$ are disjoint 
subsets of $V(G)$, 
we say $A$ is {\em complete} to $B$ if every vertex in $A$ is complete to $B$,
and $A$ is {\em anticomplete} to $B$ if there are no edges between $A$ and $B$.
The result just proved will be used in combination with the following:

\begin{thm}\label{pinkverts}
	Let $G$ be a graph, let $A,B$ be disjoint subsets of $V(G)$, and let $k,\ell, t\ge 1$ be integers. 
	Suppose that 
	\begin{itemize}
		\item for each $T\subseteq A$ with $|T|=t$, the set of vertices in $B$ complete to $T$ has 
			chromatic number at most $\ell$;
		\item every vertex in $B$ has at least $t^{k-1}$ neighbours in $A$; and
		\item $\chi(B)>k |A|^{2k-1}\ell$.
	\end{itemize}
Then there exist distinct $a_1\LL a_k\in A$, and disjoint subsets $L_1\LL L_k\subseteq B$,
each with chromatic number more than $\ell$, such that for all $i,j\in \{1\LL k\}$ and all $v\in L_j$,
$a_i$ is adjacent to $v$ if and only if $i=j$.
\end{thm}
\Proof
We proceed by induction on $k$. Suppose first that $k=1$.
Since each vertex in $B$ has a neighbour in $A$, and $\chi(B)>|A|\ell$, there is a vertex $a_1\in A$ such that the set $L_1$ of neighbours of
$a_1$ in $B$ has chromatic number more than $\ell$, and so the theorem holds. Thus we may assume that $k\ge 2$ and the
theorem holds for $k-1$. Since $B\ne \emptyset$, and each vertex in $B$ has at least $t^{k-1}$ neighbours in $A$, 
it follows
that $|A|\ge t^{k-1}$.

For each $v\in B$,
let $N_v$ be the set of neighbours of $v$ in $A$. Since $|N_v|>0$, there are at most $|A|$
possibilities for $|N_v|$; and so there exist $C\subseteq B$ with 
$$\chi(C)\ge \chi(B)/|A|>k|A|^{2k-2}\ell$$ 
such that the sets $N_v(v\in C)$ all have the same cardinality.
For each $v\in C$, let $S_v$ be the set of all $u\in C$ with $|N_u\setminus N_v|\ge t^{k-2}$. 
\\
\\
(1) {\em $\chi(S_v)> (k-1) |A|^{2k-2}\ell$ for each $v\in C$.}
\\
\\
Since $|N_v|\ge t^{k-1}$, there exist pairwise disjoint subsets $T_1\LL T_{t^{k-2}}$ of $N_v$, each of 
cardinality $t$. If $u\in C\setminus S_v$, then since $N_u, N_v$ have the same cardinality,
it follows that $|N_v\setminus N_u|< t^{k-2}$, and so $u$ is complete to one of the sets $T_1\LL T_{t^{k-2}}$.
For $1\le i\le t^{k-2}$, the set of vertices in $C$ complete to $T_i$ has chromatic number
at most $\ell$ by hypothesis; and since $C\setminus S_v$ is the union of $t^{k-2}$
of these sets, it follows that $\chi(C\setminus S_v)\le t^{k-2}\ell$. Since $\chi(C)> k |A|^{2k-2}\ell$, 
and $|A|\ge t^{k-1}\ge t$,
we deduce that
$$\chi(S_v)> k |A|^{2k-2}\ell-t^{k-2}\ell\ge (k-1) |A|^{2k-2}\ell.$$ 
This proves (1).

\bigskip

Let $a_1\LL a_k\in A$ be distinct. A {\em feathering} for the sequence $(a_1\LL a_k)$ is a sequence $(L_2\LL L_k)$ of pairwise
disjoint subsets of $C$, each with chromatic number more than $\ell$, 
such that for all $i\in \{1\LL k\}$, all $j\in \{2\LL k\}$,  and all $u\in B_j$,
$a_i$ is adjacent to $u$ if and only if $i=j$.
For each $v\in B$, let us say a {\em tail} for $v$ is a sequence $(a_1\LL a_k)$ of distinct vertices in $A$, such that
$v$ is adjacent to $a_1$ and nonadjacent to $a_2\LL a_k$, and there is a feathering for $(a_1\LL a_k)$.
\\
\\
(2) {\em For each $v\in B$, there is a tail for $v$.}
\\
\\
Since every $u\in S_v$ has a non-neighbour in $N_v$, and $\chi(S_v)>(k-1) |A|^{2k-2}\ell$, there exists $a_1\in N_v$
such that the set $B'$ of  vertices in $S_v$ nonadjacent to $a_1$  has chromatic number more than 
$(k-1) |A|^{2k-3}\ell$.
From the inductive hypothesis on $k$, applied with $A,B,k$ replaced by $A\setminus N_v$, $B'$, $k-1$ respectively, 
there exist distinct $a_2\LL a_k\in A\setminus N_v$, and disjoint subsets $L_2\LL L_k\subseteq B'$,
each with chromatic number more than $\ell$, such that for all $i\in \{2\LL k\}$,
$a_i$ is complete to $L_i$ and anticomplete to all the other $L_j$. Since $a_1$ has no neighbours 
in $B'$, it follows that $(a_1\LL a_k)$ is a tail for $v$. This proves (2).

\bigskip

For every sequence $(a_1\LL a_k)$ of distinct vertices in $A$, let $M(a_1\LL a_k)$ be the set of $v\in C$ such that
$(a_1\LL a_k)$ is a tail for $v$. By (2), $C$ is the union of the sets $M(a_1\LL a_k)$  over all choices of 
$(a_1\LL a_k)$;
and since there are at most 
$|A|^k$ such sequences $(a_1\LL a_k)$, it follows that $\chi(M(a_1\LL a_k))\ge \chi(C)|A|^{-k}$
for some choice of $(a_1\LL a_k)$. Let $L_1=M(a_1\LL a_k)$; then $\chi(L_1)>\ell$, since
$\chi(C)|A|^{-k}> (k-1) |A|^{k-2}\ell\ge \ell$. Let 
$(L_2\LL L_k)$ be a feathering for $a_1\LL a_k$. (The latter exists since $M(a_1\LL a_k)\ne \emptyset$
and so $a_1\LL a_k$ is a tail for some $v\in C$.) Since every vertex in $L_1$ is adjacent to $a_1$, and the
vertices in $L_2\LL L_k$ are all nonadjacent to $a_1$, it follows that $L_1\LL L_k$ are pairwise disjoint.
Then $a_1\LL a_k$ and $L_1\LL L_k$ satisfy the theorem. 
This proves \ref{pinkverts}.~\bbox


\section{Cores and their neighbourhoods}

For each integer $s\ge 1$, let $H_s$ be the tree with $1+s+s^2$ vertices, in which some vertex has degree $s$ and
all its neighbours have degree $s+1$.  Every tree of radius two
is an induced subgraph of $H_s$ for some choice of $s\ge 2$, and so to prove \ref{mainthm} in general, it suffices to 
prove it in the case that $H=H_s$ for some $s\ge 2$.
Thus, we need to show that for every integer $s\ge 2$ and 
every integer $d\ge 2$, there is a polynomial $f$ such that $\chi(G)\le f(\tau_d(G))$
for every $H_s$-free graph $G$. We might as well assume that $f$ is {\em increasing}, that is, $f(y)\ge f(x)$
for all $y\ge x\ge 0$, and {\em integral}, that is, all its coefficients are integers. 
(It is tempting to assume further that $f$ is of the form $f(x)=x^c$ for some $c$, but we cannot, because $\tau_d(G)$ might
be zero or one.)

To prove \ref{mainthm}, 
we will proceed by induction on $d$, with $s$ fixed. When $d=1$, the result is trivial, and when $d=2$ it follows
from \ref{poly1}, so we assume that $d\ge 2 $ and the result holds for $d$, and we will prove that it holds for $d+1$.
In summary, then, we have:
\begin{thm}\label{hyps}
	{\bf Hypothesis A:} \begin{itemize}
		\item $s,d\ge 2$ are integers. 
		\item $f$ is an increasing integral polynomial such that $\chi(G)\le f(\tau_{d}(G))$ for every 
			$H_s$-free graph $G$, and $f(t)\ge t$ for all integers $t\ge 1$. 
	\end{itemize}
\end{thm}
(We may assume the statement $f(t)\ge t$ without loss of generality, and it will be convenient later.)
We must show that there is a polynomial (and therefore there is an increasing integral polynomial) $f'$ such 
        that if $G$ is $H_s$-free then
$\chi(G)\le f'(\tau_{d+1}(G))$; or equivalently, that there is a polynomial $f''$ such that
if $G$ is $H_s$-free and
$t>\tau_{d+1}(G)$ is an integer, then $\chi(G)\le f''(t)$. Thus, to complete the proof of \ref{mainthm}, we will prove:

\begin{thm}\label{summary0}
	Assuming Hypothesis A, there is a polynomial $f_1$ such 
if $G$ is $H_s$-free and
$t>\tau_{d+1}(G)$ is an integer, then $\chi(G)\le f_1(t)$.
\end{thm}

We need the following, a special case of theorem 3.2 of~\cite{poly1}:
\begin{thm}\label{poly1better}
	Let $s\ge 2$ be an integer. If $G$ is an $H_s$-free graph, and $t>\tau_2(G)$ is an integer, then 
	$$\chi(G)\le (s(s^2+s+1)t)^{120(s^2+s+1)}.$$
\end{thm}

For integers $w,d\ge 1$,  a
subgraph $C$ of $G$ is a {\em $(w,d)$-core} if $V(C)$ is the
disjoint union of $d$ sets that are stable in $C$ (but not necessarily stable in $G$), each of
cardinality $w$ and pairwise complete in $C$. We call these
sets the {\em parts} of the core. A core is {\em stable} if each of its parts is stable in $G$.

\begin{thm}\label{core}
        Assuming Hypothesis A, let $G$ be an $H_s$-free graph, let $t>\tau_{d+1}(G)$ be an integer and let $w\ge 1$
        be an integer.
	If $G'$ is an induced subgraph of $G$ with 
	$$\chi(G')> f\left(\left(s(s^2+s+1)t\right)^{120(s^2+s+1)}w\right)$$ 
	then $G'$ has a stable $(w,d)$-core.
\end{thm}
\Proof Let $T=\left(s(s^2+s+1)t\right)^{120(s^2+s+1)}w$. Since $\chi(G')>f(T)$, it follows that $G'$ has a $(T,d)$-core $C$. 
Let $A_1\LL A_d$ be the parts of $C$.
For $1\le i\le d$, since $t>\tau_{d+1}(G)$, and the parts of $C$ all have
cardinality at least $t$, it follows that $\tau_2(G[A_i])<t$, and so by \ref{poly1better}, 
$$\chi(A_i)\le (s(s^2+s+1)t)^{120(s^2+s+1)}.$$
Hence there is a stable subset $B_i$ of $A_i$ of cardinality $w$, since $T\ge (s(s^2+s+1)t)^{120(s^2+s+1)}w$.
But then $B_1\cupcup B_d$ induces a stable $(w,d)$-core. This proves \ref{core}.~\bbox

\begin{thm}\label{colourdense}
	Assuming Hypothesis A, 
	let $G$ be an $H_s$-free graph, let $w,d\ge 1$ be integers, and let $t>\tau_{d+1}(G)$ be an integer.
	Let $C$ be a stable $(w,d)$-core, and let $P$ be the set of all vertices in $G$
	that have at least $t^{s-1}$ neighbours in some part of $C$ and have a non-neighbour in $V(C)$. 
	(Thus $V(C)\subseteq P$.) Then 
	$$\chi(P)\le sd^2w^{2s}(f(t)+ 2^{s^{2d+2}}t^{ds+s^2+s}).$$
\end{thm}
\Proof
Define
$q(t)=f(t)+ 2^{s^{2d+2}}t^{ds+s^2+s}$.
Now let $G, C,P,t$ be as in the theorem. Every vertex in $P$ has at least $t^{s-1}$ neighbours in some part of $C$,
and a non-neighbour in some part of $C$, and we may choose these two parts to be different. Since there
are only $d^2w$ choices of a part of $C$ and a vertex in a different part, 
there is a part $A$ of $C$, and a vertex $v_0$ of a different part of $C$, such that $\chi(B)\ge \chi(P)/(d^2w)$, where
$B$ is the set of vertices in $P$ that are nonadjacent to $v_0$ and have at least $t^{s-1}$ neighbours in $A$.
(Thus $A\cap B=\emptyset$, since $A$ is stable.) For each $T\subseteq A$
with $|T|=t$, the set of vertices in $B$ that are complete to $T$ contains no $(t,d)$-core (since $G$ contains no 
$(t,d+1)$-core), and so has chromatic number at most $f(t)\le q(t)$. 

Suppose for a contradiction that 
$$\chi(B)>s|A|^{2s-1}q(t).$$ 
Then by \ref{pinkverts}, taking $\ell=q(t)$, 
there exist distinct $a_1\LL a_k\in A$, and disjoint subsets $L_1\LL L_k\subseteq B$,
each with chromatic number more than $q(t)$ and hence with cardinality at least $2^{s^{2d+2}}t^{ds+s^2+s}$, such that for all 
$i,j\in \{1\LL k\}$ and all $v\in L_j$,
$a_i$ is adjacent to $v$ if and only if $i=j$.
By \ref{levels2},
there exist $I\subseteq \{1\LL k\}$ with $|I|=s$, and a subset $X_i\subseteq L_i$ for each $i\in I$, where
$\bigcup_{i\in I}X_i$ is a stable set, and $|X_i|= s$ for each $i\in I$.
But then the subgraph induced on $\{v_0\}\cup \bigcup_{i\in I}(\{a_i\}\cup X_i)$ is isomorphic to  $H_s$, a contradiction.

This proves that $\chi(B)\le s|A|^{2s-1}q(t)$, and since $|A|=w$ and $\chi(B)\ge \chi(P)/(d^2w)$, it follows that
$\chi(P)\le  sd^2w^{2s}q(t)$. This proves \ref{colourdense}.~\bbox


\section{Templates}

Our proof of \ref{mainthm} follows the ``template'' approach used by Kierstead and Penrice in~\cite{kierstead},
but modified to make the numbers polynomial.
Assuming Hypothesis A, let $w,t\ge 1$ be integers, and let us define a {\em $(w,t)$-template} in a graph $G$ to be 
a pair $(C,P)$ such that 
\begin{itemize}
	\item $C$ is a stable $(w,d)$-core in $G$;
	\item $P\subseteq V(G)$ with $V(C)\subseteq P$; 
	\item for every vertex $v\in P\setminus V(C)$ there is a part $A$ of $C$ such that $v$ has at least $st^{s-1}$ neighbours 
		in $A$; and
	\item for every vertex $v\in P\setminus V(C)$ there is a part $A$ of $C$ such that $v$ has at least $\lfloor w/t\rfloor$ non-neighbours in $A$.
\end{itemize}
A {\em $(w,t)$-template sequence}
in $G$ is a sequence $(C_i,P_i)\;(i\in I)$ of $(w,t)$-templates, where $I$ is a set of integers, such that
the sets $P_i\; (i\in I)$ are pairwise disjoint, and 
for all $i,j\in I$ with $i<j$, every vertex of $P_j$ either has fewer than $st^{s-1}$ neighbours in each part of $C_i$, or
has fewer than $\lfloor w/t\rfloor$ non-neighbours in each part of $C_i$ (that is, adding this vertex to $P_i$ would violate the
definition of a template).

The method of proof is, we will let $w$ be some appropriately large polynomial, and assume that $G$ is $H_s$-free 
and $t>\tau_{d+1}(G)$; and greedily choose a $(w=w(t),d)$-template sequence  $(C_i,P_i)\;(1\le i\le n)$ in $G$, 
where each $P_i$ is as large as possible among the set of vertices that have so far not been used, and with 
$n$ maximum. It follows that the set of vertices not in any of the templates of the sequence has bounded chromatic 
number, and so it remains to bound the chromatic number of the union of the templates. Each template has bounded
chromatic number by \ref{colourdense}, but we need to control the edges between templates in the sequence, to
bound the chromatic number of their union. Let us assume that $G$ has very large chromatic number. If we partition
the template sequence into a bounded number of other sequences, one of them will induce a subgraph that still has
very large (not quite so large) chromatic number. By this process we can make successively nicer template sequences,
still inducing large chromatic number, until eventually we will obtain a contradiction (we will obtain 
a template sequence in which for each template, each of its vertices has only a bounded number of neighbours 
in other templates of the sequence.) 

If $(C_i,P_i)\;(i\in I)$ is a $(w,t)$-template sequence in a graph $G$, its {\em vertex set} is 
$\bigcup_{i\in I} P_i$; its {\em support} is the subgraph of $G$ induced on its vertex set;
and its {\em chromatic number} is the chromatic
number of its support.

We will often use the following lemma:
\begin{thm}\label{dicolour}
	Let $D$ be a directed graph in which every vertex has out-degree at most $d$. Then $V(D)$ can be partitioned into $2d+1$ stable sets.
\end{thm}
\Proof
Every non-null subgraph $H$ has at most $d|H|$ edges, and so has a vertex $v$ such that the sum of its indegree 
and outdegree is at most $2d$. Hence the undirected graph underlying $D$ has degeneracy at most $2d$, and so is 
$(2d+1)$-colourable. This proves \ref{dicolour}.~\bbox

Let us say a $(w,t)$-template sequence $(C_i,P_i)\;(i\in I)$ with vertex set $U$ in a graph $G$ is 
\begin{itemize}
	\item {\bf 1-nice} 
if 
for each $i\in I$, there is no vertex $v\in U\setminus V(C_i)$ that has fewer than $\lfloor w/t\rfloor$ non-neighbours in 
each part of $C_i$.
\item {\bf 2-nice} if it is 1-nice
and for each $v\in U$, there are fewer than $s$ values of $i\in I$ such that $v\notin V(C_i)$
and $v$ has at least $s^3t^{s-1}$ neighbours in some part of $C_i$.
\item {\bf 3-nice} if it is 2-nice
and for all distinct $i,j\in I$, every vertex in $C_i$ has fewer than
$s^3t^{s-1}$ neighbours in each part of $C_j$.
\item {\bf 4-nice} if it is 3-nice
and for each $v\in U$, there are fewer than $(dt)^s$ values of $i\in I$
such that $v\notin P_i$ and $v$ has a neighbour in $V(C_i)$.
\item {\bf 5-nice} if it is 4-nice
and for all distinct $i,j\in I$, there are no edges between $V(C_i)$ and $V(C_j)$.
\item {\bf 6-nice} if it is 5-nice
and for all $v\in U$, there are fewer than $2(dt)^{2s}+(dt)^s$ values of $i\in I$
such that $v\notin P_i$ and $v$ has a neighbour in $P_i$.
\item {\bf 7-nice} if it is 6-nice
and for all $i\in I$, there are fewer than $60dws(dt)^{5s}$ values of $j\in I\setminus \{i\}$ such that some $v\in P_i$
has at least $(dt)^s$ neighbours in $P_j$.
\item {\bf 8-nice} if it is 7-nice and for all $i\in I$ and $v\in P_i$, $v$ has fewer than $3(dt)^{3s}$ neighbours in $U\setminus P_i$.
\end{itemize}
We will show that if $G$ has large chromatic number, then so does some 1-nice template sequence. Then 
for $i = 2\LL 8$ in turn, we will deduce
that some $i$-nice template sequence has large chromatic number; and finally, we will show by applying \ref{colourdense}
that when $i=8$, this is impossible. It will follow that $G$ has bounded chromatic number.
We begin with:
\begin{thm}\label{1nice}
	Assuming Hypothesis A, let $G$ be an $H_s$-free graph, let $t>\tau_{d+1}(G)$ be an integer, and let
	$w\ge 1$ be an integer. 
	If every 1-nice $(w,d)$-template sequence in $G$ has chromatic number at most $k$, then 
	$$\chi(G)\le f\left(\left(s(s^2+s+1)t\right)^{120(s^2+s+1)}w\right)+2kt.$$
\end{thm}
\Proof  We observe first:
\\
\\
(1) {\em Let $C$ be a $(w,d)$-core in $G$, and let $Q$ be the set of vertices $v\in V(G)\setminus V(C)$
that have fewer than $\lfloor w/t\rfloor$ non-neighbours in each part of $C$. Then $|Q|< t.$}
\\
\\
Suppose not; then there is a set $T$ of $t$ vertices of $G$, with $T\cap V(C)=\emptyset$, such that every vertex in $T$
has fewer than $\lfloor w/t\rfloor$ non-neighbours in each part of $C$. Let $A$ be a part of $C$;
then at most $(w/t-1)t$ vertices of $A$ have a non-neighbour in $T$, and since $|A|=w$, there are
$t$ vertices
in $A$ that are complete to $T$. Since this is true for each part of $C$, it follows that $G$ contains a $(t,d+1)$-core,
contradicting that $t>\tau_{d+1}(G)$. This proves (1).

\bigskip

Choose an integer $n\ge 0$, maximum such that there is a sequence $C_1,P_1,C_2,P_2\LL C_n, P_n$ 
with the following properties:
\begin{itemize}
	\item $C_1, C_2\LL C_n$ are stable $(w,d)$-cores of $G$, and $P_1\LL P_n$ are subsets of $V(G)$;
	\item the sets $P_1, P_2\LL  P_n$ are pairwise disjoint, and $V(C_i)\subseteq P_i$ for $1\le i\le n$;
	\item for $1\le i\le n$, $P_i$ consists of $V(C_i)$ together with all vertices $v$ of $G\setminus (P_1\cupcup P_{i-1})$ such that 
			$v$ has at least $st^{s-1}$ neighbours
                in some part of $C_i$ and 
        $v$ has at least $\lfloor w/t\rfloor$ non-neighbours in some part of $C_i$.
\end{itemize}
It follows that 
$(C_i,P_i)\;(i\in \{1\LL n\})$ is a $(w,t)$-template sequence in $G$.
Let its vertex set be $U$. From the maximality of $n$, there is no stable $(w,d)$-core in $G\setminus U$, and so
$$\chi(G\setminus U) \le f\left(\left(s(s^2+s+1)t\right)^{120(s^2+s+1)}w\right)$$ 
by \ref{core}. 
Let $D$ be the digraph with vertex set $\{1\LL n\}$, in which for distinct $i,j\in \{1\LL n\}$,
$j$ is adjacent from $i$ if some vertex of $P_j$ has fewer than $\lfloor w/t\rfloor$ non-neighbours in 
each part of $C_i$. By (1), every vertex of $D$ has outdegree at most $t-1$, and so by \ref{dicolour}, $V(D)$
can be partitioned into $2t$ stable sets $I_1\LL I_{2t}$. For $1\le j\le 2t$, $(C_i,P_i)\;(i\in I_j)$
is a $(w,t)$-template sequence in $G$, and it is 1-nice from the definition of $D$. Consequently each
such template sequence has chromatic number at most $k$, from the hypothesis; and so $\chi(U)\le 2kt$.
Hence 
$$\chi(G)\le \chi(U)+\chi(G\setminus U)\le 2kt+ f\left(\left(s(s^2+s+1)t\right)^{120(s^2+s+1)}w\right).$$
This proves \ref{1nice}.~\bbox

We observe:
\begin{thm}\label{1niceprop}
Assuming Hypothesis A, let $G$ be an $H_s$-free graph, let $t>\tau_{d+1}(G)$ be an integer, and let $w\ge 1$ be an integer.
Let $(C_i,P_i)\;(i\in I)$ be a 1-nice $(w,t)$-template sequence in a graph $G$. Then 
	\begin{itemize}
		\item for all $i,j\in I$
			with $i<j$, every vertex of $P_j$ has fewer than $st^{s-1}$ neighbours in each part of $C_i$.
		\item for all distinct $i,j\in I$, every vertex of $P_j$ has at least $\lfloor w/t\rfloor$
	non-neighbours in some part of $C_i$.
	\end{itemize}
\end{thm}
\Proof
Let $i,j\in I$ be distinct, and let $v\in P_j$. 
Since the sequence is 1-nice, $v$ has at least $\lfloor w/t\rfloor$ non-neighbours in some part of $C_i$;
and if $i<j$, then 
from the definition of a template sequence, 
$v$ has fewer than $st^{s-1}$ neighbours in each part of $C_i$.
This proves \ref{1niceprop}.~\bbox

\bigskip

\begin{thm}\label{2nice}
	Assuming Hypothesis A, let $G$ be an $H_s$-free graph, let $t>\tau_{d+1}(G)$ be an integer, and let 
	$w\ge (s+1)t+s^3t^{s}$ be an integer.
        Then every 1-nice $(w,d)$-template sequence in $G$ is 2-nice.
\end{thm}

\Proof Let $(C_i,P_i)\;(i\in I)$ be a 1-nice $(w,d)$-template sequence in $G$ that is not 2-nice, with vertex set $U$, 
suppose that $v\in U$ and $I'\subseteq I$ such that $|I'|=s$, and for each $i\in I'$, $v\notin V(C_i)$ and 
$v$ has at least $s^3t^{s-1}$ neighbours in some part of $C_i$.

For each $i\in I'$, there is a part $A_i$ of $C_i$ such that
$v$ has at least $s^3t^{s-1}$ neighbours in $A_i$, and a part $B_i$ of $C_i$ such that $v$ has at least
$\lfloor w/t\rfloor$ non-neighbours in $B_i$; and since $C_i$ has at least two parts, all with the same cardinality,
we may choose $A_i, B_i$ distinct.

We define $a_i\in A_i$ and $Y_i\subseteq B_i$ with $|Y_i|=s$ for $i\in I'$ inductively as follows. 
Assume that $i\in I'$, and 
$a_j$ and $Y_j$ are defined for all $j\in I'$ with $j>i$. 
Let $X=\bigcup_{j\in I',\; j>i} \{a_j\}\cup Y_j$.
Thus $|X|\le (s-1)(s+1)$. For all $j\in I'$ with $j>i$, every vertex in $V(C_j)$ has fewer than
$st^{s-1}$ neighbours in $A_i$, and since $v$ has at least $s^3t^{s-1}>(s-1)s(s+1)t^{s-1}$ neighbours in $A_i$,
there is a neighbour $a_i$ of $v$
in $A_i$ that is nonadjacent to every vertex in $X$. Similarly, since $v$ has at least
$\lfloor w/t\rfloor\ge s+ (s-1)s(s+1)t^{s-1}$ non-neighbours in $B_i$, there is a set $Y_i\subseteq B_i$ of $s$ 
vertices each nonadjacent to $v$
and each with no neighbours in $X$.
This completes the inductive definition. But then the subgraph induced on $v$ together with all the sets
$\{a_i\}\cup Y_i\;(i\in I')$ is isomorphic to $H_s$, a contradiction. This proves \ref{2nice}.~\bbox

\begin{thm}\label{3nice}
        Assuming Hypothesis A, let $G$ be an $H_s$-free graph, let $t>\tau_{d+1}(G)$ be an integer, and let
        $w\ge 1$ be an integer.
        If every 3-nice $(w,d)$-template sequence in $G$ has chromatic number at most $k$, then every 2-nice
	$(w,d)$-template sequence in $G$ has chromatic number at most $2sdwk$.
\end{thm}
\Proof
Let $(C_i,P_i)\;(i\in I)$ be a 2-nice $(w,d)$-template sequence in $G$. Let $D$ be the digraph with vertex set $I$, in
which for distinct $i,j\in I$, there is an edge of $D$ from $i$ to $j$ if some vertex of $V(C_i)$ has at least 
$s^3t^{s-1}$ neighbours in some part of $C_j$. Since $|C_i|=dw$, and the sequence is 2-nice, it follows that
every vertex of $D$ has outdegree at most $dw(s-1)$, and so $V(D)$ can be partitioned into $2dw(s-1)+1\le 2sdw$
stable sets. Each gives a 3-nice $(w,d)$-template sequence in $G$, and therefore has chromatic number at most $k$,
and so $(C_i,P_i)\;(i\in I)$ has chromatic number at most $2sdwk$. This proves \ref{3nice}.~\bbox

\bigskip

\begin{thm}\label{4nice}
        Assuming Hypothesis A, let $G$ be an $H_s$-free graph, let $t>\tau_{d+1}(G)$ be an integer, and let 
	$w\ge (s-1)s^2t^{s-1}+s^4t^{s-1}+ s$ be an integer.
        Then every 3-nice $(w,d)$-template sequence in $G$ is 4-nice.
\end{thm}
\Proof
Let $(C_i,P_i)\;(i\in I)$ be a 3-nice $(w,d)$-template sequence in $G$, with vertex set $U$, let $v\in U$, and suppose that 
there exists $I'\subseteq I$ with $|I'|\ge (dt)^s$ such that for each $i\in I'$, $v\notin P_i$ and $v$ has a neighbour 
in $V(C_i)$. 
Since $\omega(G)\le dt$, \ref{ramsey} implies that there exists $I''\subseteq I'$ with $|I''|=s$,
such that for each $i\in I'$,
$v$ has a neighbour $a_{i}\in V(C_{i})$, where the vertices $a_i\;(i\in I'')$ are pairwise nonadjacent. 
For $1\le i\le s$ let $B_i$ be a part of $C_i$ that does not contain $a_i$.
Inductively for each $i\in I''$
we choose $Y_i\subseteq B_i$ of cardinality $s$ as follows. Assume that $i\in I''$ and $Y_j$
has been defined for all $j\in I''$ with $j>i$. Let $X=\bigcup_{j\in I'', j>i}Y_j$.
Thus $|X|\le (s-1)s$.
Each vertex in $X$ has at most $st^{s-1}$ neighbours in $B_i$, and $v$ and each vertex $a_j\;(j\in I''\setminus \{i\})$
has at most $s^3t^{s-1}$ neighbours in $B_i$. Since $w-(s-1)s^2t^{s-1}-s^4t^{s-1}\ge s$, there exist a set $Y_i$ of $s$ 
distinct
vertices in $B_i$ that are nonadjacent to every vertex in $X$, and nonadjacent to 
$v$, and nonadjacent
to each vertex $a_j\;(j\in I''\setminus \{i\})$. This completes the inductive definition. But then the subgraph
induced on $\{v\}\cup \bigcup_{i\in I''}(\{a_i\}\cup Y_i)$
is isomorphic to $H_s$, a contradiction. This proves \ref{4nice}.~\bbox

\begin{thm}\label{5nice}
        Assuming Hypothesis A, let $G$ be an $H_s$-free graph, let $t>\tau_{d+1}(G)$ be an integer, and let
        $w\ge 1$ be an integer.
        If every 5-nice $(w,d)$-template sequence in $G$ has chromatic number at most $k$, then every 4-nice
	$(w,d)$-template sequence in $G$ has chromatic number at most $2dw(dt)^sk$.
\end{thm}
\Proof
Let $(C_i,P_i)\;(i\in I)$ be a 4-nice $(w,d)$-template sequence in $G$, with vertex set $U$, and let $D$
be the digraph with vertex set $I$, in which for distinct $i,j\in I$, $j$ is adjacent from $i$ if some vertex in $V(C_i)$ 
has a neighbour in $V(C_j)$. Since $|C_i|=wd$ and the sequence is 4-nice, it follows that $D$ has maximum outdegree 
less than $wd(dt)^s$. By \ref{dicolour}, $V(D)$ is the union of $2wd(dt)^s$ stable sets, each forming a 5-nice
sequence, and therefore with chromatic number at most $k$. This proves \ref{5nice}.~\bbox

\begin{thm}\label{6nice}
	 Assuming Hypothesis A, let $G$ be an $H_s$-free graph, let $t>\tau_{d+1}(G)$ be an integer, and let
        $w\ge 1$ be an integer. Then every 5-nice $(w,d)$-template sequence in $G$ is 6-nice.
\end{thm}
\Proof
Let $(C_i,P_i)\;(i\in I)$ be a 5-nice $(w,d)$-template sequence in $G$, with vertex set $U$, and let $v\in U$,
and suppose that there are at least $2(dt)^{2s}+(dt)^s$ values of $i\in I$ such that $v\notin P_i$ and $v$ has a 
neighbour in 
$P_i$. Since the sequence is 4-nice, there are fewer than $(dt)^s$ values of $i\in I$
such that $v\notin P_i$ and $v$ has a neighbour in $V(C_i)$; so there exists $I_1\subseteq I$
with $|I_1|=2(dt)^{2s}$
such that for each $i\in I_1$, $v\notin P_i$, and $v$ has a neighbour $a_i$ in $P_i$, and $v$ has no neighbour in $V(C_i)$. 
Let $D$ be the digraph with vertex set $I_1$, in which for 
distinct $i,j\in I_1$, $j$ is adjacent from $i$ if $a_{i}$ has a neighbour in $V(C_{j})$. Since the sequence
is 4-nice, $D$ has maximum outdegree at most $(dt)^s-1$, and so by \ref{dicolour}, $I_1$ can be partitioned into 
$2(dt)^s$ sets that are stable in $D$.
One of these stable sets has cardinality at least $(dt)^s$, since $|I_1|=2(dt)^{2s}$;
so there exists exists $I_2\subseteq I_1$ with $|I_2|\ge (dt)^s$, such that for all $i\in I_2$,
$v\notin P_i$, and $v$ has a neighbour $a_i$ in $P_i$, and $v$ has no neighbour in $V(C_i)$;
and for all distinct $i,j\in I_2$, $a_i$ has no neighbour in $V(C_j)$. By \ref{ramsey}, there exists $I_3\subseteq I_2$
with $|I_3|=s$, such that the vertices $a_i\;(i\in I_3)$ are pairwise nonadjacent. For each $i\in I_3$, 
since $a_i\in P_i$, there is a part of $C_i$ such that $a_i$ has at least $st^{s-1}\ge s$ neighbours in this part,
and so there exist $b^1_i\LL b^s_i\in V(C_i)$, distinct, pairwise nonadjacent, and all adjacent to $a_i$. But then the subgraph 
induced on $\{v\}\cup \bigcup_{i\in I_3}\{a_i, b^1_i\LL b^s_i\}$ is isomorphic to $H_s$, a contradiction. This proves 
\ref{6nice}.~\bbox

\begin{thm}\label{7nice}
         Assuming Hypothesis A, let $G$ be an $H_s$-free graph, let $t>\tau_{d+1}(G)$ be an integer, and let
        $w\ge 1$ be an integer. Then every 6-nice $(w,d)$-template sequence in $G$ is 7-nice.
\end{thm}
\Proof
Let $(C_i,P_i)\;(i\in I)$ be a 6-nice $(w,d)$-template sequence in $G$, with vertex set $U$, and suppose that
$h\in I$, and there exists $I_1\subseteq I\setminus \{h\}$ with $|I_1|\ge 60dws(dt)^{5s}$ such that for each $i\in I_1$ there exists 
$a_i\in P_h$ that has at least $(dt)^s$ neighbours in $P_i$.
Let $D$ be the digraph with vertex set $I_1$, where for distinct $i,j\in I_1$, $j$ is adjacent from $i$
if $a_i$ has a neighbour in $P_j$. Since the sequence is 6-nice, $D$ has maximum outdegree
less than $2(dt)^{2s}+(dt)^s\le 3(dt)^{2s}$, and so by \ref{dicolour}, there is a subset $I_2\subseteq I_1$ 
with $|I_2|\ge  |I_1|/(6(dt)^{2s})\ge  10dws(dt)^{3s}$, such that for all distinct $i,j\in I_2$, 
$a_i$ has no neighbour in $P_j$.  Since for each $i\in I_2$, $a_i$ has a neighbour in $P_i$, 
it follows that the vertices $a_i\;(i\in I_2)$ are all distinct. Since $|C_h|=dw$, there exists 
$I_3\subseteq I_2$ with $|I_3|\ge |I_2|-dw\ge 9dws(dt)^{3s}$ such that $a_i\in P_h\setminus V(C_h)$ for each $i\in I_3$. For each $i\in I_3$,
$a_i$ has a neighbour in $V(C_h)$, and since $|C_h|=dw$, there exists $I_4\subseteq I_3$ and $c\in V(C_h)$
such that 
$|I_4|\ge |I_3|/(dw)\ge 9s(dt)^{3s}$
and $a_i$ is adjacent to $c$ for each $i\in I_4$. 

There are at most $3(dt)^{2s}$ values of $i\in I_4$
such that $c$ has a neighbour in $V(C_i)\cup P_i$, so there exists $I_5\subseteq I_4$ with 
$|I_5|=|I_4|-3(dt)^{2s}\ge 6s(dt)^{3s}$ such that $c$ has no neighbour in $V(C_i)\cup P_i$ for each $i\in I_5$. 
For each $i\in I_5$,
since $a_i$ has $(dt)^s$ neighbours in $P_i$, \ref{ramsey} implies that there is a stable subset
$Y_i\subseteq P_i$ with $|Y_i|=s$ such that $a_i$ is complete to $Y_i$.
Let $D'$ be the digraph 
with vertex set $I_5$, in which for distinct $i,j\in I_5$, $j$ is adjacent from $i$ if some vertex of $Y_i$ has a 
neighbour in $Y_j$. Since each vertex in $Y_i$ has a neighbour in $P_j$
for fewer than $3(dt)^{2s}$ values of $j$, it follows that $D'$ has maximum outdegree less than
$3s(dt)^{2s}$, and so by \ref{dicolour}, there exists $I_6\subseteq I_5$
with $|I_6|\ge |I_5|/(6s(dt)^{2s})\ge (dt)^s$ such that for all distinct $i,j\in S_6$, there are no edges between 
$Y_i$ and $Y_j$. Since $|I_6|\ge (dt)^s$, \ref{ramsey} implies that there exists $I_7\subseteq I_6$
with $|I_7|=s$ such that the set of $a_i$ with
$i\in I_7$ is stable in $G$. But then the subgraph induced on 
$\{c\}\cup \bigcup_{i\in I_7} (\{a_i\}\cup Y_i)$
is isomorphic to $H_s$, a contradiction. This proves \ref{7nice}.~\bbox

\begin{thm}\label{8nice}
        Assuming Hypothesis A, let $G$ be an $H_s$-free graph, let $t>\tau_{d+1}(G)$ be an integer, and let
        $w\ge 1$ be an integer.
        If every 8-nice $(w,d)$-template sequence in $G$ has chromatic number at most $k$, then every 7-nice
        $(w,d)$-template sequence in $G$ has chromatic number at most $120dws(dt)^{5s}k$.
\end{thm}
\Proof Let $(C_i,P_i)\;(i\in I)$ be a 7-nice $(w,d)$-template sequence in $G$, with vertex set $U$. Let $D$ be the digraph 
with vertex set $I$, in which for distinct $i,j\in I$, $j$ is adjacent from $i$ if some vertex in $P_i$ has more than
$(dt)^s$ neighbours in $P_j$. Since the sequence is 7-nice, $D$ has maximum outdegree less than $60dws(dt)^{5s}$,
and so by \ref{dicolour}, $I$ is the union of $120dws(dt)^{5s}$ subsets each stable in $D$. We observe that if $I_1\subseteq I$
is stable in $D$, then for each $i\in I_1$ and each $v\in P_i$, there is no $j\in I_1\setminus \{i\}$ such that
$v$ has at least $(dt)^s$ neighbours in $P_j$; and so, since $v$ has neighbours in $P_j$
for at most $2(dt)^{2s}+(dt)^s\le 3(dt)^{2s}$ values of $j\ne i$, it follows that $v$ has fewer 
than $3(dt)^{3s}$ neighbours in $U\setminus P_i$, and so the sequence $(C_i,P_i)\;(i\in I_1)$ is 8-nice.
This proves \ref{8nice}.~\bbox

To finish this chain of reductions, we have:
\begin{thm}\label{8nicebound}
        Assuming Hypothesis A, let 
	$G$ be an $H_s$-free graph, let $t>\tau_{d+1}(G)$ be an integer, and let
	$w\ge 1$ be an integer. Then every  8-nice $(w,d)$-template sequence in $G$ has chromatic number at most 
	$$3sd^{3s+2}w^{2s-1}t^{3s}(f(t)+ 2^{s^{2d+2}}t^{ds+s^2+s}).$$
\end{thm}
\Proof
Let $(C_i,P_i)\;(i\in I)$ be an 8-nice $(w,d)$-template sequence in $G$, with vertex set $U$.  Let $G_1, G_2$ be the subgraphs 
of $G$, both with vertex set $U$, where $G_1$ is the union of the subgraphs $G[P_i]\;(1\le i\le n)$ and $G_2$
contains precisely those edges of $G[U]$ that do not belong to $G_1$. By \ref{colourdense}, 
$$\chi(G_1)\le sd^2w^{2s}(f(t)+ 2^{s^{2d+2}}t^{ds+s^2+s}),$$
and $\chi(G_2)\le 3(dt)^{3s}$ since $G_2$ has maximum degree less than $3(dt)^{3s}$ (because the sequence is 8-nice).
Taking the product colouring shows that 
$$\chi(U)\le 3sd^{3s+2}w^{2s-1}t^{3s}(f(t)+ 2^{s^{2d+2}}t^{ds+s^2+s}).$$ This proves \ref{8nicebound}.~\bbox

\bigskip

Now we can complete the proof of \ref{mainthm}, by proving \ref{summary0}, which we restate:
\begin{thm}\label{summary}
	Assuming Hypothesis A, there is a polynomial $f_1$ with the following property. Let $G$ be an $H_s$-free graph, 
	and let $t>\tau_{d+1}(G)$ be an integer.
	Then $\chi(G)\le f_1(t)$.
\end{thm}
\Proof Define 
\begin{align*}
	w(t)&=s^4t^s+s\\
	f_8(t)&=3sd^{3s+2}w^{2s-1}t^{3s}(f(t)+ 2^{s^{2d+2}}t^{ds+s^2+s})\\
	f_5(t)=f_6(t)=	f_7(t) &= 120sd^{5s+1}wt^{5s}f_8(t)\\
	f_3(t)=f_4(t)&= 2d^{s+1}wt^sf_5(t)\\
	f_2(t)&=2sdwf_3(t)\\
	f_1(t)&=f\left(\left(s(s^2+s+1)t\right)^{120(s^2+s+1)}w\right)+2tf_2(t).
\end{align*}
Thus, $w,f_8, f_7\LL f_1$ are all polynomials in $t$, since $d,s$ are constants and $f$ is a polynomial by Hypothesis A.
By \ref{8nicebound}, every  8-nice $(w,d)$-template 
sequence in $G$ has chromatic number at most $f_8(t)$.
For $i = 7,6\LL 1$ in turn, it follows that every $i$-nice $(w,d)$-template
sequence in $G$ has chromatic number at most $f_i(t)$, by applying \ref{8nice}, \ref{7nice}, \ref{6nice}, \ref{5nice},
\ref{4nice}, \ref{3nice}, \ref{2nice} respectively.
By \ref{1nice}, $\chi(G)\le f_1(t)$. 
This proves \ref{summary}.~\bbox

\section*{Acknowledgement}
Our thanks to Sophie Spirkl, who worked with us on some parts of this paper. And we would like to thank the referees for a 
thorough and very helpful report.

\end{document}